\documentclass[a4paper,12pt]{article}
\usepackage{fullpage}
\usepackage[english]{babel}
\usepackage[latin1]{inputenc}
\usepackage[T1]{fontenc}
\usepackage{float}
\usepackage{amsmath}
\usepackage{amssymb}
\usepackage{graphicx}
\usepackage{dsfont}% Pour avoir le 1 caract\'{e}ristique
\newtheorem{rque}{Remark}[section]
\newtheorem{lem}{Lemma}[section]
\newtheorem{thm}{Theorem}[section]

\newcommand{\R}{\mathbb{R}}
\newcommand{\N}{\mathbb{N}}
\newcommand{\Z}{\mathbb{Z}}

\newcommand{\U}{\overrightarrow{U}}
\newcommand{\ii}{\overrightarrow{i}}
\newcommand{\ji}{\overrightarrow{j}}
\newcommand{\ki}{\overrightarrow{k}}
\newcommand{\fig}{ F{\sc{ig}}. }
\newcommand{\vecdeux}[2]{\left(\begin{array}{c} #1\\#2 \end{array}\right)}
\newcommand{\vectrois}[3]{\left(\begin{array}{c}#1\\#2\\#3\end{array}\right)}

\newcommand{\dive}{\textbf{div}}

\newbox\hautbox \setbox\hautbox=\hbox{\vphantom{\rule[0.2cm]{0cm}{0.2cm}}} % Permet de redefinir la hauteur de ligne

\title{A kinetic scheme for pressurised flows in non  uniform closed water pipes}

\author{C. Bourdarias$^{1}$\thanks{ {\it Christian.Bourdarias@univ-savoie.fr}}, 
M. Ersoy$^1$\thanks{{\it Mehmet.Ersoy@univ-savoie.fr}} and 
S. Gerbi$^{1}$\thanks{{\it Stephane.Gerbi@univ-savoie.fr}} \\[1mm]
{\small $^{1}$ Universit\'{e} de Savoie,Laboratoire de Math\'{e}matiques,}{\small 73376 Le Bourget-du-Lac, France.}}

\date{\today}

\begin{document}
\maketitle
% \tableofcontents
\abstract{The aim of this paper is to present a kinetic numerical scheme for the computations of transient
pressurised flows in closed water
pipe with non uniform sections. Firstly, we detail the derivation of the mathematical model in curvilinear coordinates  and  we performe  a  formal asymptotic analysis. The obtained system  is  written as a conservative hyperbolic partial
differential system of equations. We obtain a kinetic interpretation of this system and we build the corresponding kinetic scheme based on an
upwinding of the source terms written as the  gradient of a  ``pseudo altitude''.
The validation is lastly performed in the case of a water hammer in an uniform pipe: we compare the numerical results provided by an industrial code used at EDF-CIH (France), which
solves the Allievi equations 
(the commonly used equation for pressurised flows in pipe) by the method of characteristics, with
those of the kinetic scheme. To validate the contracting or expanding cases, we compare the presented technique to the equivalent pipe method in the case of an immediate flow shut down in a quasi-frictionless cone-shaped pipe.}


\vspace{0.5cm}
\section{Introduction}

The presented work takes place in a more general project: the modelization of unsteady mixed flows in any kind of closed domain taking into account  the cavitation problem and air entrapment. We are interested in flows occuring in closed pipe of non uniform sections, where some parts of the flow can be free surface (it means that only a part of the pipe is filled) and other parts are pressurised (it means that the pipe is full-filled). The transition phenomenon, between the two types of flows, occurs in many situation such as storm sewers, waste or supply pipes in hydroelectric installation. It can be induced by sudden change in the boundary conditions as failure pumping. During this process, the pressure can  reach severe values and  cause damages. 

The classical Shallow Water equations are commonly used to describe  free surface flows in open channel. They are also used in the study of mixed flows  using the Preissman slot artefact (see for example \cite{CSZ97,SWB98}). However, this technic does not take into account  depressurisation phenomenon which occurs during a water hammer. We can also cite  the Allievi equations which are commonly used to describe pressurised flows. Nonetheless, the non conservative form is not well adapted to a natural coupling with the Shallow Water equations (contrary to the one presented in \cite{BEG08}).

The model for the unsteady mixed water flows in closed water pipes and a finite volume
discretization has been previously studied by two of the authors \cite{BG07} and a kinetic formulation  has
been proposed in \cite{BGG08}. This paper tends to extend naturally the work in \cite{BGG08} in the case of closed pipes with non uniform sections.

We establish, in Section \ref{model}, the model for pressurised flows in curvilinear coordinates and recall some classical properties of this model. Rewritting the source terms due to both topography and geometry into a single one that we called \emph{pseudo-altitude} term, we get a model close to the presented  one by the authors in \cite{PS01}. In Section \ref{kineticscheme},  we  present the kinetic formulation of this model that will be useful to show the main properties of the numerical scheme.  The last part is devoted to the construction of the kinetic scheme: the upwinding of
the source term due to the pseudo topography 
is performed in a close manner described by  Perthame \textit{et al.} \cite{PS01} using an energetic balance
at microscopic level. We have used  the generalized characteristics method to extend the works in \cite{BGG08} to  the kinetic scheme with \emph{pseudo-reflections}.

Finally,  we present in Section \ref{numerics} a numerical validation of this study in the uniform case by the
comparison between the resolution of this model and the resolution of the Allievi equation solved by
the industrial code \verb+belier+ used at Center in Hydraulics Engineering of Electricit\'{e} De
France (EDF) \cite{W93} for the case of critical water hammer tests. The validation in  non uniform pipes is performed in the case of an immediate flow shut down in a quasi-frictionless cone-shaped pipe. The results are compared to the equivalent pipe method \cite{A03}.

\section{Formal Derivation of the model}\label{model}
The presented model is derived from the 3D compressible Euler system  written in curvilinear coordinates, then integrated over sections orthogonal to the main flow axis (see below). We  neglect the second and third equation of the conservation of the momentum and we  get an unidirectionnal model. Then, an asymptotic analysis is performed to get a model close to the Shallow Water model (to a future coupling for the study of unsteady mixed flows \cite{BEG08}).
\subsection{The Euler system in curvilinear coordinates}
The 3D Euler system in the cartesian coordinates is written as follows 
\begin{equation}\label{Euler3D_mass_conservation}
\partial_t \rho + \dive{(\rho \U)} = 0,
\end{equation}
\begin{equation}\label{Euler3D_momentum_conservation}
\partial_t (\rho\U) + \dive{(\rho \U\otimes \U)} + \nabla p = F,
\end{equation}

\noindent where $\U(t,x,y,z)$ and $\rho(t,x,y,z))$ denotes the velocity with components $(u,v,w)$ and  the density respectively. $p(t,x,y,z)$ is the scalar pressure and $F$ the exterior strenght of gravity. 

\begin{figure}[H]
\centering{
\includegraphics[scale=0.6]{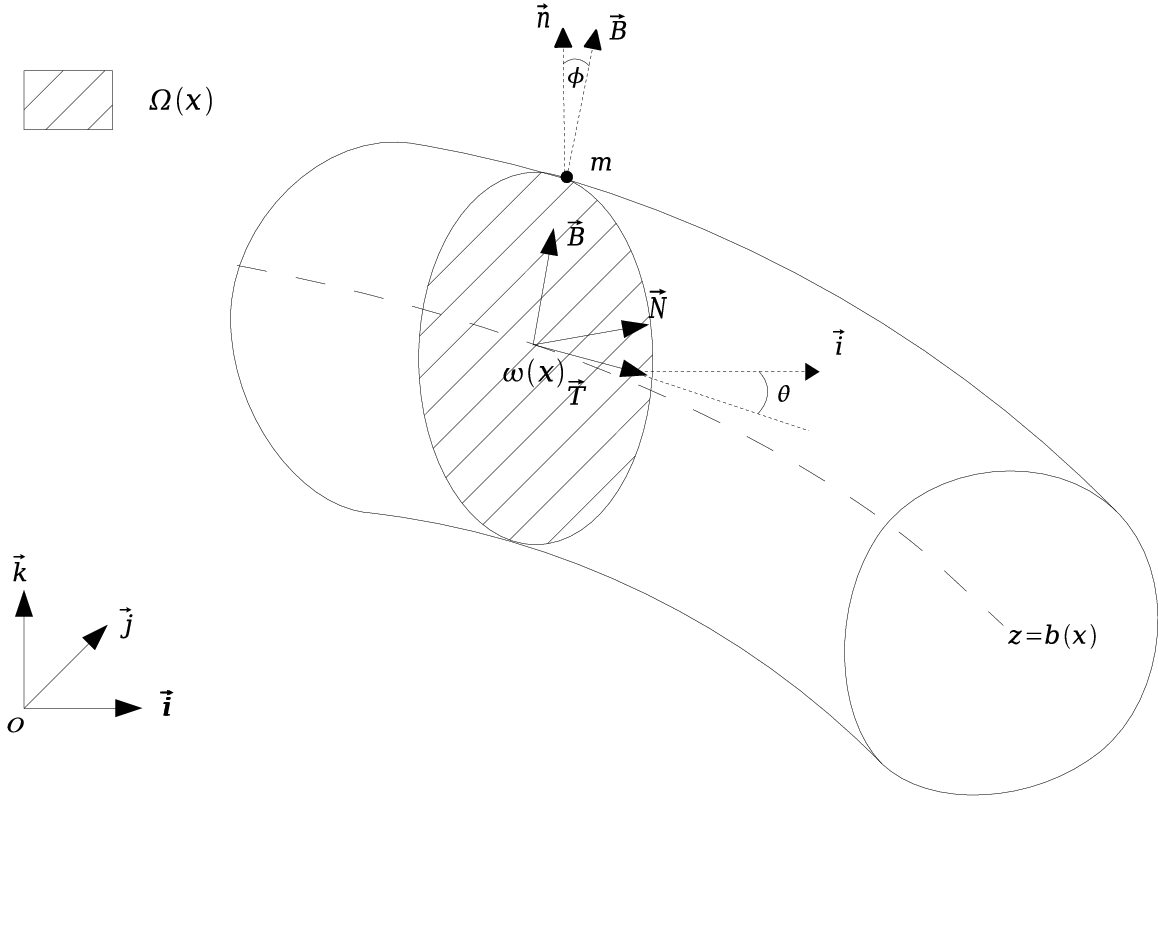}
}
\caption{Geometric characteristics of the pipe}\label{pressurised_closed_pipe_domain}
\end{figure}

\noindent We define the domain $\Omega_F$  of the flow as the union of sections $\Omega(x)$ (assumed to be simply connected compact sets) orthogonal  to some plane curve with parametrization $(x,0,b(x))$ in a convenient cartesian reference frame $(O,\ii,\ji,\ki)$ where $\ki$ follows the vertical direction; $b(x)$ is then the elevation of the point $\omega (x,0,b(x))$ over the plane $(O,\ii,\ji)$ (see \fig\ref{pressurised_closed_pipe_domain}). The curve may be, for instance, the axis spanned by the center of mass of each orthogonal section $\Omega(x)$ to the main mean flow axis, especially in the case of a piecewise cone-shaped pipe. Notice that we consider only the case of rigid pipe: the sections are only $x$-dependent.

\noindent To see the local effect induced by the geometry due to the change of sections and/or  slope, we write the  3D compressible Euler system in the  curvilinear coordinates. To this end, let us introduce the curvilinear variable defined by \\$\displaystyle X = \int_{x_0}^x \sqrt{1+(b'(\xi))^2} d\xi$ where  $x_0$ is  an arbitrary abscissa. We set $y=Y$ and we denote by $Z$ the altitude of any fluid particle $M$ in the Serret-Frenet basis  $(\overrightarrow T,\overrightarrow N,\overrightarrow B)$ at point $\omega (x,0,b(x))$:  $\overrightarrow T$ is the tangent vector, $\overrightarrow N$ the normal vector and $\overrightarrow B$ the binormal vector (see\fig\ref{pressurised_closed_pipe_domain}). Then we perform the following transformation $\mathcal{T}: (x,y,z)\rightarrow (X,Y,Z)$ and we use the following lemma (whose proof can be found in \cite{BFML07}):
\begin{lem}\label{lemma_chain_rule} Let $(x,y,z) \mapsto \mathcal{T}(x,y,z)$ be a transformation and $\mathcal{A}^{-1} = D_{(x,y,z)} \mathcal T$ the jacobian matrix of the transformation with determinant $J$.

Then, for any vector field $\Phi$ one has,
$$J D_{(X,Y,Z)} \Phi = D_{(x,y,z)}(J\mathcal{A}\Phi)
$$
In particular, for any scalar function $f$, one has
$$D_{(X,Y,Z)}f = \mathcal{A}^t D_{(x,y,z)}f
$$
\end{lem}

\noindent Let $(U,V,W)^t$ be the components of the velocity vector  in the $(X,Y,Z)$ coordinates  in such a way  that the flow is orthogonal to the sections $\Omega(x)$. Let $R$ be the matrix defined by $ R=\left( 
\begin{array}{ccc}
 \cos\theta & 0 & \sin\theta \\
 0 & 1 & 0 \\
 -\sin\theta& 0 & \cos\theta \\
\end{array} \right)$ then: $$\left(\begin{array}{c} U\\V\\W
\end{array}\right) = R\left(\begin{array}{c} u\\v\\w
\end{array}\right).$$

\noindent Applying Lemma  \ref{lemma_chain_rule} to the mass conservation equation, we get
$$J(\partial_t \rho + \dive(\rho \U)) =0$$
$$\Longleftrightarrow$$
\begin{equation}\label{mass_conservation_curvilinear}
\partial_t (J\rho) + \partial_X(\rho U) + \partial_Y(\rho JV)+ \partial_Z(\rho JW) = 0
\end{equation}

\noindent where \begin{align} J = & \det \left(
\begin{array}{ccc}
 \left(1-Z \displaystyle \frac{d}{dX}\theta\right)\cos\theta & 0 & \sin\theta \\
 0 & 1 & 0 \\
 \left(1-Z \displaystyle \frac{d}{dX}\theta\right)\sin\theta& 0 & \cos\theta \\
\end{array}
\right).\end{align}

\noindent To get the unidirectionnal Shallow Water-like equations, we suppose that the mean flow follows the $X$-axis. Hence, we neglect the second and third equation for the conservation of the momentum. Therefore, we only perform the curvilinear transformation for the first conservation equation. To this end, multiplying the conservation of the momentum equation of System (\ref{Euler3D_momentum_conservation}) by $J\vectrois{\cos\theta}{0}{\sin\theta}$ and using Lemma \ref{lemma_chain_rule}  yields:
$$
J\vectrois{\cos\theta}{0}{\sin\theta}\left(\partial_t( \rho\U) + \dive(\rho \U\otimes\U) + \nabla p = -\rho\nabla(\overrightarrow g.\overrightarrow{OM})\right).
$$
\noindent It may be rewritten as:
\begin{equation}\label{momentum_conservation_curvilinear}
\begin{array}{c}
\partial_t(J\rho U) + \partial_X(\rho U^2) + \partial_Y(\rho JUV^2)+ \partial_Z(\rho JUW) + \partial_X p \\ 
= -\rho Jg\sin\theta+\rho U W  \displaystyle \frac{d}{dX}(\cos \theta)
\end{array}
\end{equation} where $\overrightarrow{OM}$ denotes the position of any particule $M$ in the local basis \newline
\noindent $(\overrightarrow T,\overrightarrow N,\overrightarrow B) \textrm{ at point } \omega (x,0,b(x)).$

\noindent Finally, in the $(X,Y,Z)$ coordinates the system reads: 
\begin{equation}\left\{
\begin{array}{rcl}
\partial_t (J\rho) + \partial_X(\rho U) + \partial_Y(\rho JV)+ \partial_Z(\rho JW) &=& 0\\ 
& & \label{Euler_system_curvilinear}\\
\partial_t(J\rho U) + \partial_X(\rho U^2) + \partial_Y(\rho JUV^2)+ \partial_Z(\rho JUW)+ \partial_X p \\ = -\rho Jg\sin\theta  +\rho U W  \displaystyle \frac{d}{dX}(\cos \theta)\end{array}\right.\end{equation}

\begin{rque}\rm ~ 
Notice that $\kappa(X)=\displaystyle\frac{d}{dX}\theta$ is the algebric curvature of the axis at $\omega(x)$ and the function $J(X,Y,Z) = 1-Z \kappa(X)$  only depends on variables $X,Z$. Morever, we assume $J>0$ in $\Omega_F$ which corresponds to a reasonnable geometric hypothesis.
Consequently, $J$ defines a diffeomorphism and thus the performed transformation is admissible.
\end{rque}

\noindent We recall that the main objective is to obtain a formulation close to the Shallow Water equation in order to couple the two models in a natural way (in a close manner described in \cite{BG07}). The direct integration of Equations (\ref{Euler_system_curvilinear}) over $\Omega(x)$ gives a model which is not useful, due to the term $J$, to perform a natural coupling with the Shallow Water model \cite{BEG08} for non uniform pipes. Setting $\epsilon = H/L $ a small parameter (where $H$ and $L$ are two characteristics dimensions along $\ki$ and $\ii$ axis respectively), we get  $J = 1+O(\epsilon)$. We also assume that the characteristic dimension along the $\ji$ axis is the same as $\ki$.
We introduce the others characteristics dimensions $T,P,\overline U,\overline V,\overline W$ for time, pressure and velocity repectively and the dimensionless quantities as follows: $$\widetilde U = U/ \overline U,\, \widetilde V = \epsilon V/ \overline U,\, \widetilde W = \epsilon {W}/{\overline U},\, $$ $$\widetilde X =
{X}/{L},\, \widetilde Y = {Y}/{H},\, \widetilde Z = {Z}/{H},\,
\widetilde p = {p}/{P},\widetilde\theta = {\theta},\widetilde\rho = {\rho}.$$ In the sequel, we set $P={\overline
U}^2$ and $L=T\overline U$ (i.e. we consider only laminar flow).

\noindent Under these hypothesis $J(X,Y,Z)=\widetilde J (\widetilde X,\widetilde Y,\widetilde Z)
= 1-\epsilon \widetilde Z \displaystyle \frac{d}{d\widetilde X}\theta$. So, the rescaled system~(\ref{Euler_system_curvilinear}) reads: 
\begin{equation}\label{Euler_curvilinear_rescaled}\left\{
\begin{array}{rcl}
\partial_{\widetilde t}(\widetilde J \widetilde\rho) + \partial_{\widetilde X}({\widetilde\rho\widetilde U}) +
\partial_{\widetilde Y}(\widetilde J \widetilde\rho\widetilde V)+
\partial_{\widetilde Z}(\widetilde J \widetilde\rho\widetilde W)&=&0\\
& & \\
\partial_{\widetilde t}(\widetilde J \widetilde U \widetilde\rho) + \partial_{\widetilde X}({\widetilde\rho\widetilde U}^2) +
\partial_{\widetilde Y}(\widetilde J \widetilde\rho\widetilde U \widetilde V)+
\partial_{\widetilde Z}(\widetilde J \widetilde\rho\widetilde U \widetilde W) + \partial_{\widetilde X} \widetilde p  \\ =   \epsilon \widetilde\rho\widetilde U \widetilde W \widetilde \rho(\widetilde
  X)\displaystyle-\widetilde\rho\displaystyle\frac{\sin\widetilde\theta}{{F_{r,L}}^2}
\displaystyle-  \frac{\widetilde Z \partial_{\widetilde X}(\cos\theta)}{{F_{r,H}}^2}
\end{array}\right.\end{equation}
 with
 $F_{r,M} = \displaystyle\frac{\overline U}{\sqrt{gM}}$ the Froude number along the $\ii$ axis and the $\ki$ or $\ji$ axis where $M$ is any generic variable.

\noindent Formally, when $\epsilon $ vanishes the  system
reduces to:
\begin{equation}\left\{
\begin{array}{rcl}
\partial_{\widetilde t}(\widetilde\rho) + \partial_{\widetilde X}({\widetilde\rho\widetilde U}) +
\partial_{\widetilde Y}(\widetilde\rho\widetilde V)+
\partial_{\widetilde Z}(\widetilde\rho\widetilde W)&=&0\\
& & \\
\partial_{\widetilde t}(U \widetilde\rho) + \partial_{\widetilde X}(\widetilde\rho{\widetilde U}^2) +
\partial_{\widetilde Y}(\widetilde\rho\widetilde U \widetilde V)+
\partial_{\widetilde Z}(\widetilde\rho\widetilde U \widetilde W)+ \partial_{\widetilde X} \widetilde p &=&
  -\widetilde\rho\displaystyle\frac{\sin\widetilde\theta}{{F_{r,L}}^2}\\
& &-\displaystyle
  \frac{\widetilde Z \partial_{\widetilde X}(\cos\widetilde\theta)}{{F_{r,H}}^2} \end{array}\right.\end{equation}

\noindent Finally, the system in variables $(X,Y,Z)$ that describes the slope variation and the section variation in a closed pipe reads: 
\begin{equation}\label{Euler_to_integrate}\left\{
\begin{array}{rcl}
\partial_{t}(\rho) + \partial_{X}({\rho U}) +
\partial_{Y}(\rho V)+
\partial_{Z}(\rho W)&=&0\\
& & \\
\partial_{ t}(U  \rho) + \partial_{  X}({ \rho  U}^2) +
\partial_{  Y}( \rho  U   V)+
\partial_{  Z}( \rho  U   W)+ \partial_{  X}   p &=&
  -\rho\displaystyle{g\sin \theta} \\ & & -
    Z  \displaystyle \frac{d}{dX}(g\cos \theta)
\end{array}\right.\end{equation}

\begin{rque}\rm
To take into account the friction, we add the source term $-\rho g S_f \overrightarrow T$ (described above) in the momentum equation.
\end{rque}

\subsection{Shallow Water-like equations in closed pipe}
In the following,  we  use the linearized pressure law $p = p_a + \displaystyle{\frac{\rho-\rho_0}{\beta \rho_0}}$ (see e.g. \cite{SWB98,WS78}) in which $\rho_0$ represents the density of the fluid at atmospheric pressure $p_a$ and $\beta$ the water compressibility coefficient  equal to $5.0\,10^{-10}\,m^2.N^{-1}$ in practice. The sonic speed is then given by $c = 1/\sqrt{\beta \rho_0}$ and thus $c \approx 1400\,m.s^{-1}$. The friction term is given by the Manning-Strickler law (see \cite{SWB98}),
$$S_f = K(S)U|U|\textrm{ with } K(S) = \frac{1}{K_s^2R_h(S)^{4/3}}$$ where $S=S(X)$ is the surface area of the section $\Omega(X)$  normal to the main pipe  axis (see \fig\ref{pressurised_closed_pipe_domain} for the notations). $K_s$ is the coefficient of roughness and $R_h(S) = S/P_m$ is the hydraulic radius where $P_m$ is the perimeter of $\Omega$. 

\noindent System (\ref{Euler_to_integrate}) is integrated over the cross-section $\Omega$. 
In the following, overlined letters represents  the averaged quantities over $\Omega$. For $m \in \partial \Omega$, $\displaystyle\overrightarrow n = \frac{\overrightarrow m}{|\overrightarrow m|}$ is the outward unit vector at the point $m$ in the $\Omega$-plane and $\overrightarrow m $ stands for the vector $\overrightarrow{\omega m}$ (as displayed on \fig\ref{pressurised_closed_pipe_domain}).

\noindent Following the work in \cite{BG07}, using the approximations $\overline{\rho U}\approx \overline{\rho}\overline{U},\,\overline{\rho U^2}\approx \overline{\rho}\overline{U}^2$ and Lebesgue integral formulas, the mass conservation equations becomes:
\begin{equation}\label{STV_mass}\partial_{t}(\overline\rho S) + \partial_{X}({\overline\rho q}) = \displaystyle \int_{\partial \Omega} \rho \left( U\partial_X \overrightarrow m - \overrightarrow V\right).\overrightarrow n\, ds,\end{equation}
where $q = S \overline U$ is the discharge of the flow and the velocity  $\overrightarrow V = (V,W)^t$ in the $(\overrightarrow N, \overrightarrow B)$-plane.

\noindent The equation of the conservation of the momentum becomes 
\begin{equation}\label{STV_momentum}
\begin{array}{lll}
\partial_{t}(\overline\rho q) + \partial_{X}(\displaystyle\frac{\overline\rho q^2}{S}+c^2\overline\rho S) &=&
-\displaystyle{g\overline\rho S \sin \theta}+c^2\overline\rho \frac{d S}{dX}\\	
 & -&   \overline\rho S\overline Z  \displaystyle \frac{d}{dX}(g\cos\theta)\\
 &  +&   \displaystyle \int_{\partial \Omega} \rho U \left(U\partial_X \overrightarrow m - \overrightarrow V\right).\overrightarrow n\, ds
\end{array}
\end{equation}

\noindent The integral terms appearing in (\ref{STV_mass}) and (\ref{STV_momentum}) vanish, as the pipe is infinitely rigid, i.e. $\Omega = \Omega(X)$ (see \cite{BG07} for the dilatable case). It follows the  non-penetration condition: $$\vectrois{U}{V}{W}.\overrightarrow N = 0.$$ 

\noindent Finally, omitting the overlined letters except $\overline{Z}$, we obtain the equations for pressurised flows under the form
\begin{equation}\label{SVLike}\left\{
\begin{array}{rcl}
\partial_{t}(\rho S) + \partial_{X}({\rho q}) &=&0\\
\partial_{t}(\rho q) + \partial_{X}(\displaystyle\frac{\rho q^2}{S}+c^2\rho S) &=&
  -\displaystyle{\rho S g\sin \theta }-
    \rho S \overline Z  \displaystyle \frac{d}{dX}(g\cos \theta) + c^2 \rho \frac{d S}{dX}
\end{array}\right.\end{equation}
where the quantity $\overline Z$ is the $Z$  coordinate of the center of mass. 
\begin{rque}
In the case of a circular section pipe, we choose the plane curve $(x,0,b(x))$ as the mean axis and we get obviously $\overline Z = 0$.
\end{rque}

\noindent Now, following \cite{BG07}, let us introduce the conservative variables $A = \displaystyle\frac{\rho S}{\rho_0}$ the \emph{equivalent wet area} and the \emph{equivalent discharge}  $Q =\displaystyle A U$. Then dividing System (\ref{SVLike}) by $\rho_0$ we get:
\begin{equation}\label{SVLike_M_D}\left\{
\begin{array}{rcl}
\partial_{t}(A) + \partial_{X}(Q) &=&0\\
&  & \\
\partial_{t}(Q) + \partial_{X}(\displaystyle\frac{Q^2}{A}+c^2 A) &=&
  -\displaystyle{g A \sin\theta}-
A \overline Z   \displaystyle \frac{d}{dX}(g\cos \theta) + \\ &  &\displaystyle c^2 A \frac{d}{dX}\ln(S)
\end{array}\right.\end{equation}

\begin{rque}\rm
 This choice of variables is motivated by the fact that this system is formally closed to the Shallow Water equations with topography source term in non uniform pipe. Indeed, the Shallow water equations for non uniform pipe reads \cite{BEG08}: 
$$
\left\{\begin{array}{lll}
\partial_{t}A +\partial_X Q &=& 0\\
\displaystyle\partial_t Q +\displaystyle\partial_X\left(\frac{Q^2}{ A}+\displaystyle g\cos\theta I_1 \right)&=&-gA\sin\theta -A(h-I_1(A)/A)\displaystyle \frac{d}{dX}(g\cos \theta)\\ & &  +\displaystyle g\cos\theta I_2\end{array}\right.
$$
\noindent where the terms $gI_1 \cos\theta$, $I_2 \cos\theta$, $(h-I_1(A)/A)$ are respectively the equivalent terms to $c^2 A$, $\displaystyle c^2 A \frac{d}{dX}\ln(S)$, $\overline{Z}$ in System (\ref{SVLike_M_D}). The quantities $I_1$, $I_2$, $(h-I_1(A)/A)$  denotes respectively the classical term of hydrostatic pressure, the pressure source term induced by the change of geometry and the $Z$  coordinate of the center of mass. 
\noindent Finally, the choice of these unknowns leads to a natural coupling between the pressurised and free surface model (called PFS-model presented  by the authors in \cite{BEG08}).
\end{rque}

\noindent To close this section, let us give the classical properties of System (\ref{SVLike_M_D}): 
\begin{thm}[frictionless case]\label{Classical_prop}
\begin{enumerate}
\item[] 
\item The system (\ref{SVLike_M_D}) is stricly hyperbolic for $A(t,X)>0$.
\item For smooth solutions, the mean velocity $U = Q/A$ satisfies 
\begin{equation}
\partial_t U + \partial_X \left(\displaystyle\frac{ U^2}{2} +c^2\ln (A/S) +g\Phi_{\theta}+gZ \right) = 0
\end{equation}
where $\Phi_{\theta}(X) = \displaystyle \int_{X_0}^X \overline Z(\xi)\displaystyle\frac{d}{dX}\cos\theta(\xi)\,d\xi$ for any arbitrary $x_0$ and $Z$ the altitude term defined by $\partial_X Z = sin\theta $. 
The quantity $\displaystyle \frac{ U^2}{2} +c^2\ln (A/S) +g\Phi_{\theta}+gZ$
is also called the total head. 
\item The still water steady states for $U = 0$ is given by 
\begin{equation}
c^2\ln (A/S) +g\Phi_{\theta}+gZ = 0.
\end{equation}
\item It admits a mathematical entropy $$E(A,Q) = 
\displaystyle\frac{Q^2}{2A}+c^2A\ln(A/S)+ gA\Phi_{\theta} +gAZ$$
which satisfies the entropy inequality $$\partial_t E +\partial_X \big((E+c^2 A )U\big) \leqslant 0$$
\end{enumerate}
\end{thm}

\begin{rque}\rm 
\begin{itemize}
 \item[]
 \item If we consider the friction term, we have  for smooth solutions: $$
\partial_t  U + \partial_X \left(\displaystyle\frac{U^2}{2} +c^2\ln (A/S) +g\Phi_{\theta}+gZ \right) = -g K(S) {U}|{U}|
$$ and the previous entropy equality reads
$$\partial_t E +\partial_X \big((E+c^2 A )U\big) = -gAK(S) U^2 |U| \leqslant 0$$

\item If we introduce $\widetilde Z$ the so-called \emph{pseudo altitude} source term  given by $$\widetilde Z= Z+\Phi_{\theta}-\displaystyle \frac{c^2}{g} \ln(S)$$ (where $\Phi_{\theta}$ is defined in Theorem \ref{Classical_prop}), we can rewrite System (\ref{SVLike_M_D}) in the simpler form, close to the classical Shallow Water formulation:
\begin{equation}\label{SVLike_M_D_SV}\left\{
\begin{array}{rcl}
\partial_{t}(A) + \partial_{X}({Q}) &=&0\\
\partial_{t}(Q) + \partial_{X}\left(\displaystyle\frac{Q^2}{A}+p(X,A)\right) +g\partial_X\widetilde Z&=&0
\end{array}\right.\end{equation} where $p(X,A)=c^2 A$.
\end{itemize}
\noindent This reformulation allows us to perform an analysis close to the presented one by the autors in \cite{PS01} in order to write the kinetic formulation.
\end{rque}

\section{The kinetic model}\label{kineticscheme}
We present in this section  the kinetic formulation (see e.g. \cite{P02}) for pressurised flows in water pipes modelized by System (\ref{SVLike_M_D_SV}). To this end, we introduce a smooth real function $\chi$ such that $$\chi(w)=\chi(-w)\geq 0,\,\displaystyle \int_{\R}\chi(w)\,dw=1,\,\displaystyle \int_{\R}w^2\chi(w)\,dw=1$$ and defines the Gibbs equilibrium as follows $$\mathcal M (t,x,\xi) =\displaystyle \frac{A}{c}\chi \left(\frac{\xi-U}{c}\right)$$ which represents the density of particles at time $t$, position $x$ and the kinetic speed $\xi$. Then we get the following kinetic formulation:
\begin{thm}\label{thm_Equiv_Kinet_Macro}
The couple of functions $(A,Q)$ is a strong solution of the Shallow Water-like system (\ref{SVLike_M_D_SV}) if and only if $\mathcal M$ satisfies the kinetic transport equation
\begin{equation}\label{kinetic_equation}
\partial_t \mathcal M + \xi \partial_X \mathcal M - g\partial_X\widetilde Z  \partial_{\xi} \mathcal M = K(t,x,\xi)
\end{equation}
for some collision kernel $K(t,x,\xi)$ which admits  vanishing moments up to order $1$ for a.e (t,x). 
\end{thm}

\textbf{Proof of Theorem \ref{thm_Equiv_Kinet_Macro}.} We get easily the above result since the following macro-microscopic relations holds
\begin{equation}\label{Kinetic_relation_M}
A = \displaystyle\int_{\R} \mathcal M(\xi)\,d\xi
\end{equation}

\begin{equation}\label{Kinetic_relation_D}
Q = \displaystyle\int_{\R} \xi\mathcal M(\xi)\,d\xi
\end{equation}
 
\begin{equation}\label{Kinetic_relation_Flux}
\displaystyle\frac{Q^2}{A}+c^2 A = \displaystyle\int_{\R} \xi^2 \mathcal M(\xi)\,d\xi
\end{equation}

\begin{flushright}
$\square$
\end{flushright}

\noindent The reformulation of System (\ref{SVLike_M_D}) and the above theorem has the advantage to give only one linear transport equation for $\mathcal M$ which can be easily discretised (see for instance \cite{PS01,PT90}). Morever, the following results hold: 
\begin{thm}\label{thm_minimization}
Let us consider the minimization problem $ \min \mathcal{E}(f) $  under the constraints $$ f >0,\quad \displaystyle\int_{\R} f(\xi)\,d\xi = A,\quad \displaystyle\int_{\R} \xi f(\xi)\,d\xi = Q$$ where the kinetic functional energy  is defined by  
$$\mathcal{E}(f) = \displaystyle\int_{\R}\frac{\xi^2}{2}f(\xi)+c^2f(\xi)log(f(\xi) )+c^2f(\xi)log(c\sqrt{2\pi})+g\widetilde Z f(\xi)\,d\xi.$$  
Then the minimum  is attained by the function $\mathcal M (t,x,\xi)=\displaystyle \frac{A}{c}\chi \left(\frac{\xi-U}{c}\right)$ where $\chi(w) = \displaystyle \frac{1}{\sqrt{2\pi}} \exp \left(\displaystyle\frac{-w^2}{2}\right) $ a.e.  

Morever, the minimal energy is $$\mathcal{E}(\mathcal M) = E(A,Q) = \displaystyle\frac{Q^2}{2A}+c^2 A\ln A + gA\widetilde Z $$
and $\mathcal M$ satisfies the still water steady state equation for $ U = 0$, that is, $$\xi\partial_X \mathcal M -g\partial_X \widetilde Z \partial_{\xi}\mathcal M = 0.$$ 
\end{thm}

\textbf{Proof of Theorem \ref{thm_minimization}} One may easily verify that  $f=\mathcal M$ is a solution of the minimization problem. Under the hypothesis $f >0$ the functionnal $\mathcal{E}(f)$ is strictly convex which ensures the unicity of the minimum. Furthermore, by a direct computation, one has $\mathcal{E}(\mathcal M) = E$.

The minimum $\mathcal M$ of the functionnal $\mathcal{E}(f)$ satisfies the still water steady state for $U=0$,$$\xi\partial_X \mathcal M -g\partial_X \widetilde Z \partial_{\xi}\mathcal M = 0.$$
Since $\partial_X \mathcal M = \displaystyle \frac{\partial_X A}{c}\chi\left(\displaystyle\frac{\xi}{c}\right)$, $\partial_{\xi} \mathcal M =  \displaystyle \frac{A}{c^2}\chi'\left(\displaystyle\frac{\xi}{c}\right)$, denoting $w = \xi/c$, we get $$w \partial_X A \chi(w)-g\partial_X\widetilde Z \displaystyle\frac{A}{c}\chi'(w) = 0. $$
On the other hand, the still water steady state at macroscopic level is given by $$c^2\ln (A) +g\widetilde Z=cst,$$ and so one has $g\partial_X\widetilde Z = -c^2\partial_X(\ln A)$. 
Finally, we get the following  ordinary differential equation  $$ w\chi(w)+\chi'(w)=0.$$ which gives the result.
\begin{flushright}
$\square$
\end{flushright}

\section{The kinetic scheme with pseudo-reflections}\label{section_Kinetic_scheme}
This section is devoted to the construction of the numerical kinetic scheme and its properties. The numerical scheme is obtained by  using a flux splitting method on the previous kinetic formulation (\ref{kinetic_equation}). The source term due to the pseudo topography  $\partial_X \widetilde Z$ is upwinded in a close manner described by  Perthame \textit{et al.} \cite{PS01} using an energetic balance
at the microscopic level. 
\noindent In the sequel, for the sake of simplicity, we consider the space domain infinite.

\noindent Let us consider the discretization $(m_i)_{i\in\Z}$ of the spatial domain with $$m_i = (X_{i-1/2},X_{i+1/2}),\,h_i = X_{i+1/2}-X_{i-1/2}$$ which are respectively the cell and mesh size for  $i\in\Z$. Let $\Delta t^n= t_{n+1}-t_n,\,n\in\N$ be the timestep.

\noindent Let  $\mathcal{U}_i^n =(A_i^n,Q_i^n), U_i^n = \displaystyle \frac{Q_i^n}{A_i^n}$ be respectively the approximation of the mean value of $(A,Q)$ and the velocity $U$ on $m_i$ at time $t_n$. 

\noindent Let $\mathcal M_i^n(\xi) = \displaystyle \frac{A_i^n}{c}\chi\left(\frac{\xi-U_i^n}{c}\right)$ be  the approximation of the microscopic quantities and $\widetilde{Z}_i \mathds{1}_{m_i}(X)$ be the piecewise constant representation of the pseudo-altitude $\widetilde{Z}$. Then, integrating System (\ref{SVLike_M_D_SV}) over $m_i\times[t_n,t_{n+1}]$,  we get: 
\begin{equation}\label{SolutionWithKineticScheme}
\mathcal U_i^{n+1} = \mathcal U_i^n -\frac{\Delta t^n}{h_i}\left(F_{i+1/2}^- -  F_{i-1/2}^+\right)
\end{equation}   
where 
\begin{equation}\label{MicroMacroFlux}F_{i+1/2}^{\pm} = \displaystyle\frac{1}{\Delta t^n}\int_{t_n}^{t_{n+1}}F\left(\mathcal{U}(t,X_{i+1/2}^{\pm})\, dt \right)\end{equation} are the interface fluxes with $F(A,Q) = (Q,Q^2/A +c^2A)^t$.

\noindent Now, it remains to define an approximation $F_{i\pm 1/2}^{\pm}$ of the flux at the points $X_{i\pm 1/2}$. To this end, we use the  kinetic formulation (\ref{kinetic_equation}).\newline
Assume that the discrete macroscopic vector state  $\mathcal{U}_i^n$ is known at time $t_n$.
We consider the  following problem  
\begin{equation}\label{Mf_kinetic_equation}
\left\{\begin{array}{lll}
\partial_t f + \xi\partial_X f - g \partial_X(\widetilde{Z}) \,\partial_{\xi}f =0 &\quad (t,X,\xi)\in[t_n,t_{n+1}]\times m_i\times{\R}\\
f(t_n,X,\xi) =\mathcal{M}(t_n,X,\xi) & \quad (X,\xi)\in m_i\times{\R}
\end{array}\right.
\end{equation}
where $\mathcal M(t_n,X,\xi) = \mathcal{M}_i^n(\xi)$ in the cell $m_i$.
It is discretized as follows (since it is a linear transport equation)
\begin{equation}\label{Mf_kinetic_equation_discret} 
\forall i\in \mathbb{Z}, \, \forall n \in \N,\,\quad f_i^{n+1}(\xi) =
\mathcal{M}_i^{n}(\xi) -\xi\frac{\Delta t^n}{h_i}\left\{\mathcal{M}_{i+1/2}^-(\xi)-\mathcal{M}_{i-1/2}^+(\xi) \right\}
\end{equation}
where $\mathcal{M}_{i\pm 1/2}^{\pm}$ denotes the interface density equilibrium (computed in section \ref{Interface_equilibrium_densities}). Finally, 
we set
\begin{equation}\label{discrete_macro_micro_relation}
\mathcal{U}_i^{n+1}=\displaystyle \int_{\R} \vecdeux{1}{\xi} f_i^{n+1}(\xi)\, d\xi 
\end{equation}
and
$$\mathcal M_i^{n+1} = \displaystyle \frac{M_i^{n+1}}{c}\chi\left(\frac{\xi-U_i^{n+1}}{c}\right).$$

\begin{rque}\rm
We can understand Equation (\ref{Mf_kinetic_equation}) as follows: let us consider the following problem,
\begin{equation}\label{Mf_kinetic_equation2}
\left\{\begin{array}{lll}
\partial_t f + \xi\partial_X \mathcal{M} - g \partial_X(\widetilde{Z}) \,\partial_{\xi}\mathcal{M} =0 &\quad (t,X,\xi)\in[t_n,t_{n+1}]\times m_i\times{\R}\\
f(t_n,X,\xi) =\mathcal{M}(t_n,X,\xi) & \quad (X,\xi)\in m_i\times{\R}.
\end{array}\right.
\end{equation}
Assuming that  $\mathcal M(t,X,\xi)$ is known on $[t_n,t_{n+1}]\times m_i\times{\R}$ leads to the same discretization (\ref{Mf_kinetic_equation_discret}) of Equation (\ref{Mf_kinetic_equation}). Hence the numerical scheme (\ref{Mf_kinetic_equation_discret}) avoids to compute explicitely  the collision kernel $K$ at the microscopic level. Indeed, substracting Equation (\ref{kinetic_equation}) to Equation (\ref{Mf_kinetic_equation2}), we get:
$$
\partial_{t}(\mathcal{M}-f)(\xi)=K(t,x,\xi).
$$ 
Then, integrating the previous identity in time $t$ and $\xi$ yields to:
$$
\int_{\R}\vecdeux{1}{\xi}f(\xi)\, d\xi = \mathcal{U}.
$$ In other words, using the numerical scheme (\ref{Mf_kinetic_equation_discret}) and the macroscopic-microscopic relation (\ref{discrete_macro_micro_relation}) is a manner to perform all collisions at once and to recover exactly the macroscopic unknows $(A,Q)$.
\end{rque}

\noindent Now to complete the numerical kinetic scheme, it remains to define the microscopic fluxes $\mathcal{M}_{i\pm 1/2}^{\pm}$ appearing in equation (\ref{Mf_kinetic_equation_discret}) introduced by the choice of the constant piecewise representation of the pseudo-altitude term $\widetilde Z$.

\subsection{Interface equilibrium densities} \label{Interface_equilibrium_densities}
\noindent To compute the interface equilibrium densities, we use the generalized characteristics method. Let $s \in (t_n,t_{n+1})$ be a time variable and $f$ the solution of the kinetic equation (\ref{Mf_kinetic_equation}). Let $i\in\Z$, $t\in (t_{n},t_{{n}+1})$ and $\xi_l$, $\xi_r$  be respectively the kinetic speed of a particle at time $t$ on each side of the interface $X_{i+1/2}$.
The characteristic curves $\Xi(s)$ and $X(s)$ of the kinetic transport equation (\ref{Mf_kinetic_equation}) satisfies the following equations:
\begin{equation}\label{ODE_curve}
\left\{\begin{array}{ccl}
\displaystyle\frac{d \Xi}{ds} &=& -g\partial_x \widetilde Z (X(s))\\
 & & \\
\displaystyle\frac{d X}{ds} &=& \Xi(s)
\end{array}\right.
\end{equation}
\noindent  where the final conditions are defined by 
\begin{equation}\label{ODE_curve_cond}
\left\{\begin{array}{ccl}
\Xi(t) &=& \xi\\
X(t) &=& X_{i+1/2}
\end{array}\right.
\end{equation} for some  constant $\xi$ defined later.
\noindent By a straightforward computation, we get the following mechanical conservation law:
\begin{equation}\label{loi_de_conservation}
\displaystyle\frac{d}{ds}\left(\frac{\Xi(s)^2}{2} + g \widetilde Z(s)\right) = 0.
\end{equation}
\noindent Since $\widetilde{Z}$ is a piecewise constant function,  the solution $\Xi$ of the ordinary differential equation (\ref{ODE_curve}) is a piecewise constant solution. So, we need to define an admissible jump condition to get only physical solutions of the problem (\ref{ODE_curve}). Thanks to the relation (\ref{loi_de_conservation}), we get the jump condition: 
$$ \left[\Xi^2\right]= \left[2g \widetilde Z\right]
$$
\noindent that is also: 
\begin{equation}\label{conservation_nrj}
\displaystyle\frac{\xi_l^2}{2} -\frac{\xi_r^2}{2} = g \Delta \widetilde Z_{i+1/2}
\end{equation}
\noindent where $ \Delta \widetilde Z_{i+1/2}$ is such that $$\widetilde Z_{i+1}-\widetilde Z_i= \Delta \widetilde Z_{i+1/2}\delta_{X_{i+1/2}}$$ with $\delta_a$ is the Dirac mass at  point $a$. The quantity $\Delta \widetilde Z_{i+1/2}$ is the \emph{potential bareer}.

\noindent Next, solving System (\ref{ODE_curve}) on $m_i\times (t_n,t_{n+1})$ with the final conditions :
\begin{equation}\label{ODE_curve_cond_complete}
\left\{\begin{array}{ccl}
\Xi(t) &=& \xi_l\\
X(t) &=& X_{i+1/2}
\end{array}\right.,
\end{equation}
\noindent we get 
\begin{equation}\label{ODE_solution}
\Xi(s) = \xi_l \textrm{ and } X(s) = \xi_l(s-t_{n+1})+X_{i+1/2}.
\end{equation}

\begin{figure}[H]
\centering{
\includegraphics[scale=0.65]{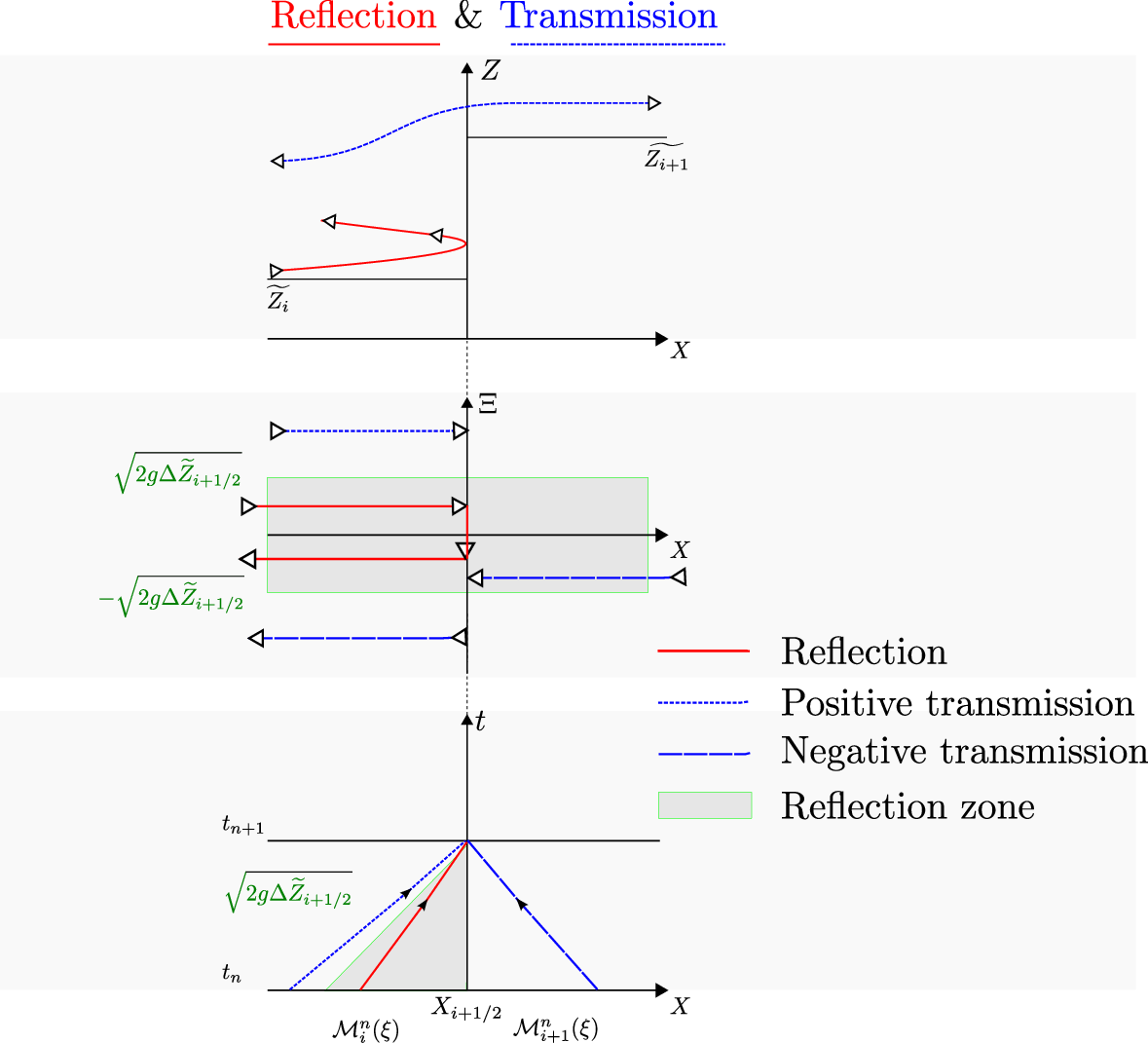}
}
\caption{The potential bareer: transmission and reflection of particle}\label{potential_bareer}
\small{Top: the physical configuration}\\
\small{Middle: the characteristic solution in $(X,\Xi)$-plane}\\
\small{Bottom: the characteristic solution in $(X,t)$-plane}
\end{figure}

\noindent Due to the jump condition (\ref{conservation_nrj}) and the sign of the kinetic speed,  we distinguish three admissible cases as displayed on \fig{\ref{potential_bareer}}.

\begin{itemize}
 \item[-] \noindent The case $\xi_l>0$ corresponds to the positive transmission (this means that the particle comes from the left) and we deduce from Equalities (\ref{ODE_solution}) that the left microscopic flux $\mathcal M_{i+1/2}^-(\xi)$ is equal to $\displaystyle\mathcal M_i^n(\xi)$.

 \item[-] The case $\xi_l<0$ and $\xi_l^2- 2{g \Delta \widetilde Z_{i+1/2}}<0$ is  the so-called  \emph{reflection} case. The condition $\xi_l^2- 2{g \Delta \widetilde Z_{i+1/2}}<0$ says simply that the slope $\xi_l$ of the $X$ solution (\ref{ODE_solution}) cannot exceed $\sqrt{2g \Delta \widetilde Z_{i+1/2}}$ (as displayed on \fig{\ref{potential_bareer}} (bottom)) and so the flux $\mathcal M_{i+1/2}^-(\xi)$ is given by $\mathcal M_i^n(-\xi)$. Physically, since the particle with the kinetic speed $\xi_l$, under the previous kinetic condition, has not enough energy to overpass the bareer, it is reflected with the kinetic speed $-\xi_l$.

\item[-] The last case is when $\xi_l<0$ and $\xi_l^2- 2{g \Delta \widetilde Z_{i+1/2}}>0$. This case corresponds to the negative transmission: this means we take into account the particles coming from the right side with negative kinetic speed. Contrary to the reflection case, the constraint on the $X$ slope is limited by $\xi_l>-\sqrt{2g \Delta \widetilde Z_{i+1/2}}$ and we get as solution 
$\displaystyle \mathcal M_{i+1}^n\left(\displaystyle-\sqrt{\xi^2-2g\Delta\widetilde Z_{i+1/2}}\right)$. From a physical point of view, the observed particle at the left of the interface comes from the right side with a kinetic speed $\xi_r<0$ where $\displaystyle \xi_r = -\sqrt{\xi_l^2-2g\Delta\widetilde Z_{i+1/2}}$, taking into account the gain or loss of potential energy through the bareer (as displayed on \fig{\ref{potential_bareer}} (bottom)).

\end{itemize}

\noindent Finally, adding the previous results we obtain:
\begin{equation}\label{interface densities}
\begin{array}{lll}
\mathcal M_{i+1/2}^{-}(\xi) &=&  \overbrace{\displaystyle \mathds{1}_{\xi>0}\mathcal M_i^n(\xi)}^{{positive \,\,\, transmission}}+\overbrace{\mathds{1}_{\xi<0,\xi^2-2g\Delta \widetilde Z_{i+1/2}<0}\mathcal M_i^n(-\xi)}^{reflection}\\ &+& \underbrace{\displaystyle \mathds{1}_{\xi<0,\xi^2-2g\Delta \widetilde Z_{i+1/2}>0}\mathcal M_{i+1}^n\left(\displaystyle-\sqrt{\xi^2-2g\Delta\widetilde Z_{i+1/2}}\right)}_{negative \,\,\,transmission}\\
 & & \\
 & & \\
\mathcal M_{i+1/2}^{+}(\xi) &=& \overbrace{\displaystyle \mathds{1}_{\xi<0}\mathcal M_{i+1}^n(\xi)}^{negative \,\,\, transmission}+\overbrace{\mathds{1}_{\xi>0,\xi^2+2g\Delta \widetilde Z_{i+1/2}<0}\mathcal M_{i+1}^n(-\xi)}^{reflection}\\ &+& \underbrace{ \displaystyle \mathds{1}_{\xi>0,\xi^2+2g\Delta \widetilde Z_{i+1/2}>0}\mathcal M_{i}^n\left(\displaystyle\sqrt{\xi^2+2g\Delta \widetilde Z_{i+1/2}}\right)}_{positive \,\,\, transmission}
\end{array}
\end{equation}
The microscopic flux at the right of the interface is obtained following a same approach.

\subsection{Numerical properties}\label{subsection_numerical_properties}
\noindent We present some numerical properties of the macroscopic scheme (\ref{SolutionWithKineticScheme})-(\ref{MicroMacroFlux}), namely the stability and the preservation of the still water steady state. The stability of the kinetic scheme  depends on a kinetic CFL condition $$\displaystyle\frac{\Delta t^n}{\max_{i} h_i} \xi < 1 ,\,\forall \xi$$ and so, on the support of the maxwellian function (e.g. we see that from the microscopic fluxes in Subsection \ref{Interface_equilibrium_densities}). The support of the  maxwellian function computed in Theorem \ref{thm_minimization} is not compact, then  the stability condition cannot be satisfied. 
Therefore, in the sequel, we will consider the particular Gibbs equilibrium $\chi(w)= \displaystyle \frac{1}{2\sqrt{3}}\mathds{1}_{[-\sqrt{3},\sqrt{3}]}(w)$ introduced by the authors in \cite{ABP00} and used in \cite{BGG08} in the case of pressurised flows in uniform closed pipe. 

\noindent Let us present the numerical properties of the scheme (\ref{Mf_kinetic_equation})-(\ref{interface densities}),
\begin{thm}\label{Numer_prop_with_indic}
\begin{enumerate}
 \item[]
 \item  Assuming the CFL condition $$\displaystyle\frac{\Delta t^n}{\max_{i\in\Z} h_i} \max_{i\in\Z}\left(|\overline U_i^n| +\sqrt{3}c\right)< 1,$$ the numerical scheme (\ref{Mf_kinetic_equation})-(\ref{interface densities}) keeps the wet equivalent area $A$ positive.
\item The  still water steady state is preserved:  $$\displaystyle U_i^n=0,\quad \frac{c^2}{g} \ln(\rho_i^n)+\widetilde Z_i = cst$$
\end{enumerate}
\end{thm}

\noindent \textbf{Proof of Theorem \ref{Numer_prop_with_indic}.} (It is similar to the one obtained in \cite{PS01}) Let us suppose $A_i^{n}>0$ for all $i\in\Z$ and $n\in\N$.
\noindent Let $\xi_{\pm} = \max(0,\pm\xi)$ be the positive or negative part of any real and $\sigma=\displaystyle\frac{\Delta t^n}{\max_{i} h_i}$, Equation (\ref{Mf_kinetic_equation}) reads:
$$\begin{array}{lll}
 f_i^{n+1}(\xi)  & \geqslant&  (1-\sigma |\xi|)\mathcal M_i^n(\xi)\\
                 &   & + \sigma \xi_+ \left(\mathds{1}_{\xi^2+2g\Delta \widetilde Z_{i+1/2}<0}\mathcal M_i^n(-\xi)\right.\\ 
                 &   & \left.+  \mathds{1}_{\xi^2+2g\Delta \widetilde Z_{i-1/2}>0}\mathcal M_{i-1}^n\left(\displaystyle\sqrt{\xi^2+2g\Delta \widetilde Z_{i+1/2}}\right)\right)  \\
		 &   &  +  \sigma \xi_- \left(\mathds{1}_{\xi^2-2g\Delta \widetilde Z_{i+1/2}<0}\mathcal M_i^n(-\xi)\right.\\  
                 &   &  \left.+ \mathds{1}_{\xi^2-2g\Delta \widetilde Z_{i-1/2}>0}\mathcal M_{i+1}^n\left(\displaystyle-\sqrt{\xi^2-2g\Delta \widetilde Z_{i+1/2}}\right)\right) 
\end{array}$$
\noindent Since the support of the $\chi$ function is compact, we get $$f_i^{n+1}(\xi)>0 \text{ if } |\xi-u_j^n|<\sqrt{3}c,\,\forall j\in\Z$$ which implies $|\xi|<|u_j^n|+\sqrt{3}c$. Using the CFL condition $\sigma |\xi|\leq 1$, we  get the result. Morever, since $f_i^{n+1}$ is a sum of positive term, we obtain $f_i^{n+1}>0$, hence the wet equivalent area at time $t^{n+1}$ is positive, i.e. $$A_i^{n+1} = \displaystyle \int_{\R} f_i^{n+1}(\xi)\,d\xi>0.$$

\noindent To prove  the second point, we distinguish cases $\xi>0$ and $\xi<0$ to show the equality $\mathcal M_{i+1/2} = \mathcal M_{i-1/2}$. Using the jump condition (\ref{conservation_nrj}), we easily obtain $f_i^{n+1} = \mathcal M_i^{n}$  which gives the result.
 \begin{flushright}
 $\square$
\end{flushright}

\noindent Now let us also remark that 
the kinetic scheme (\ref{Mf_kinetic_equation_discret})-(\ref{interface densities}) is wet equivalent area conservative .
Indeed, let us denote the first component of the discrete fluxes 
$\left(F_A\right)^\pm_{i+1/2}$:
\begin{equation*}
\left(F_A\right)^\pm_{i+1/2}
{:=} \int_\R \xi \mathcal{M}^\pm_{i+1/2}(\xi)\,d\xi 
\end{equation*}
\noindent An easy computation, using the change of variables $w^2 = {\xi}^2 - 2g\Delta \widetilde Z_{i+1/2}$
in the interface densities formulas defining the kinetic fluxes $\mathcal{M}_{i+1/2}^\pm$, allows us
to show that:
\begin{equation*}
\left(F_A\right)^+_{i+\frac{1}{2}} = \left(F_A\right)^-_{i+\frac{1}{2}}
\end{equation*}

\section{Numerical Validation}\label{numerics} 
\noindent The validation is  performed in the case of a soft and sharp  water hammer in an uniform pipe. Then we compare the results to the ones provided by an industrial code used at EDF-CIH (France) (see \cite{W93}), which
solves the Allievi equation 
by the method of characteristics. 
\noindent The validation in  non uniform pipes is performed in the case of an immediate flow shut down in a quasi-frictionless cone-shaped pipe. The results are then compared to the equivalent pipe method \cite{A03}.

\subsection{The uniform case}
\noindent We present now numerical results of a water hammer test.
The pipe of circular cross-section of $2 \mbox{ m}^2$
and thickness $20$ cm is $2000$ m long. The altitude of the upstream end of the
pipe is $250$ m and the slope is $5^{\circ}$. The Young modulus is $23 \, 10^9 \mbox{ Pa}$ since the
pipe is supposed to be built in concrete.
\noindent The total upstream head is 300 m. The initial downstream discharge is $10 \mbox{ m}^3/\mbox{s}$
and we cut the flow in $10$ seconds for the first test case and in $5$ seconds for the other.

\noindent We present a validation of the proposed scheme by comparing numerical results of the proposed model
solved by the kinetic scheme with the ones obtained by solving Allievi equations by the method of 
characteristics with the so-called \verb+belier+ code: an industrial code used by the engineers of the 
Center in Hydraulics Engineering of Electricit\'{e} De France (EDF) \cite{W93}.

\begin{figure}[H]
\centering{
\includegraphics[scale=1]{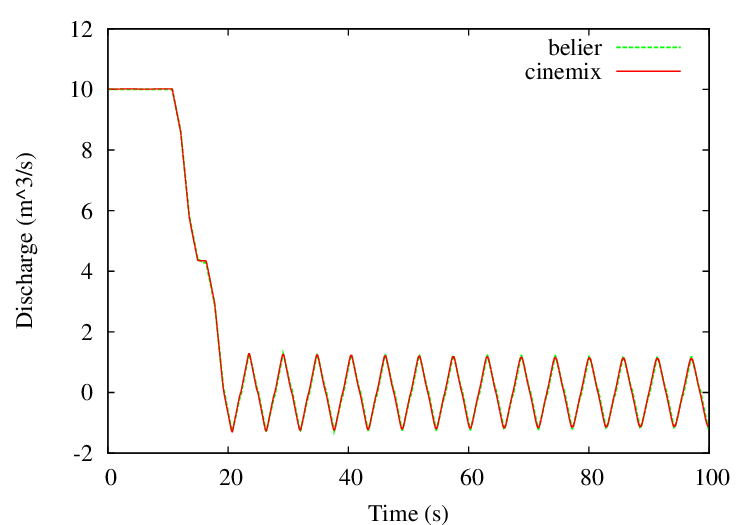}
\caption{Comparison between the kinetic scheme and the industrial code belier}
\label{order1}
{\small First case: discharge at the middle of the pipe}
}
\end{figure}

\begin{figure}[H]
\centering{
\includegraphics[scale=1]{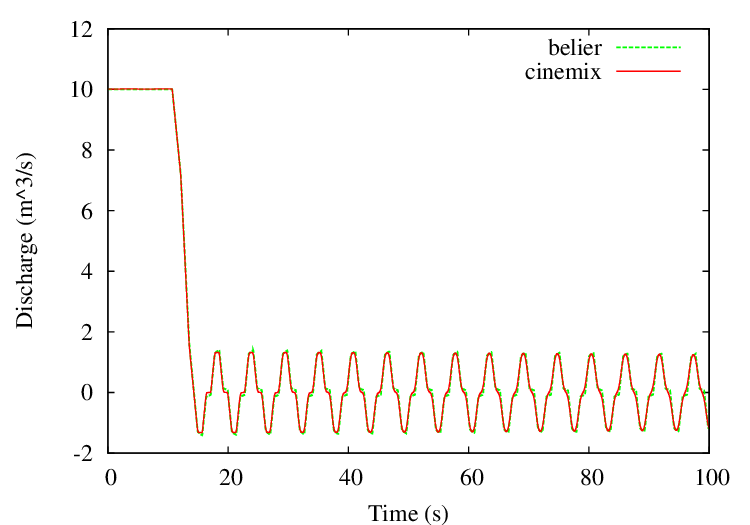}
\caption{Comparison between the kinetic scheme and the industrial code belier}
\label{order2}
{\small Second case: discharge at the middle of the pipe}
}
\end{figure}

\noindent A simulation of the water hammer test was done for a CFL coefficient equal to $0.8$ 
and a spatial discretisation of 1000 mesh points.
In the figures \fig\ref{order1} and \fig\ref{order2}, we present a comparison between the results obtained by our kinetic scheme and
the ones obtained by the ``belier'' code:
the behaviour of the discharge at the middle of the pipe. 
One can observe that the results for the proposed model are in very good agreement with the  solution of 
Allievi equations. A little smoothing effect and absorption may be probably due to the first order discretisation type.
A second order scheme may be implemented naturally and will produce a better approximation.


\subsection{The case of non uniform circular pipe}
\noindent We present a test of the proposed kinetic scheme in the case of a contracting or expanding  circular pipes  of length $L=1000\,m$. The downstream radius is kept constant, equal to $R_2 =1\,m$ and the upstream radius varies from $R_1 = 1\,m$ to $4\,m$ by steps of  $0.25\,m$. The others  paramaters are $N =300$ mesh points,  $K_S = 9000 $ (this means that the wall of the pipe is very smooth), CFL$=0.8$.
\noindent The upstream discharge before the shut-down ($1.5$ seconds) is fixed to $10\,m^3.s^{-1}$ while the upstream condition is a constant total head. 
\noindent We assume also that the pipe is rigid. Then for each value of the radius $R_1$, we compute the water hammer pressure rise at the position $x=96\, m$ of the pipe and we compare it to  the one obtained by the equivalent pipe method (see \cite{A03}).
The results are presented in \fig{\ref{equivalent_pipe}} and  show a very good agreement.

\begin{figure}[H]
\begin{center}
 \begin{tabular}{rl}
   \includegraphics[scale=0.59]{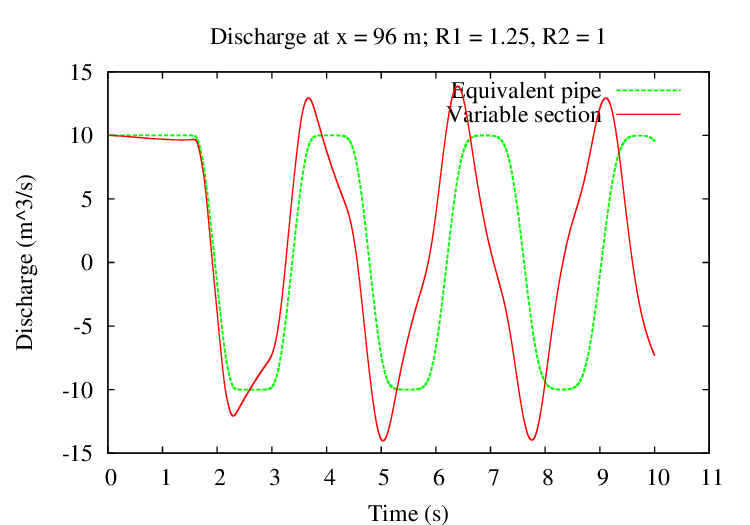} &
   \includegraphics[scale=0.59]{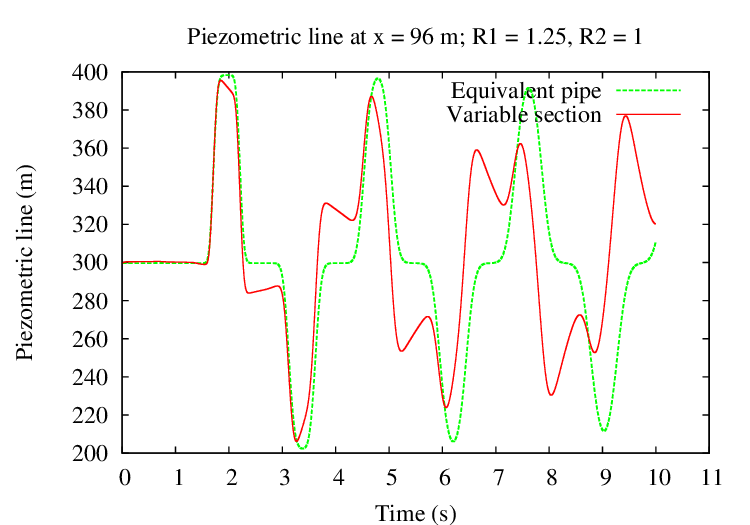} \\
 \end{tabular}
\end{center}
\caption{Discharge (left) and piezometric line (right) for $R_1 = 1.25\,m$}
\label{T12}
\end{figure}

\begin{figure}[H]
\begin{center}
 \begin{tabular}{rl}
   \includegraphics[scale=0.59]{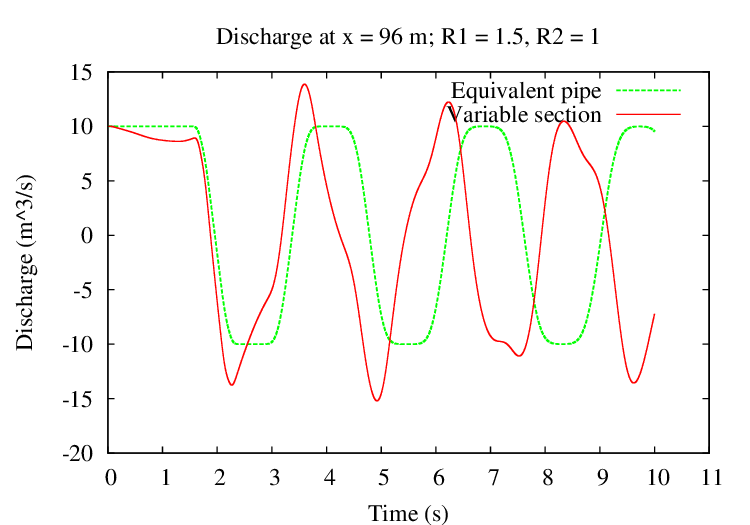} &
   \includegraphics[scale=0.59]{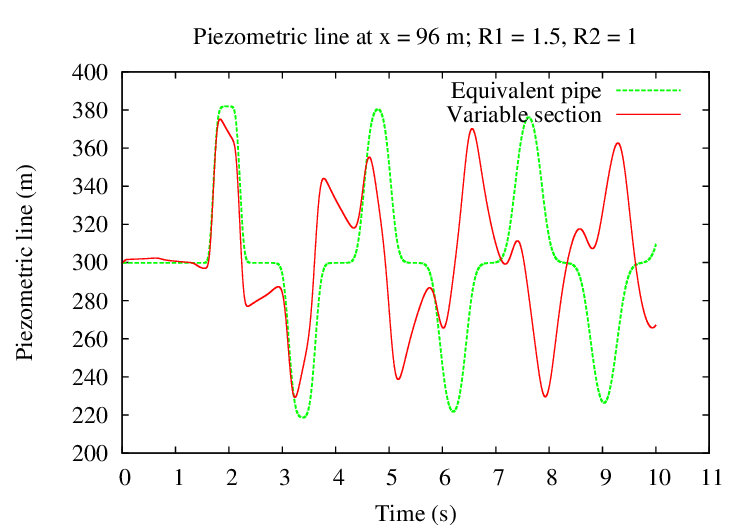} \\
 \end{tabular}
\end{center}
\caption{Discharge (left) and piezometric line (right) for $R_1 = 1.5\,m$}
\label{T11}
\end{figure}

\begin{figure}[H]
\begin{center}
 \begin{tabular}{rl}
   \includegraphics[scale=0.59]{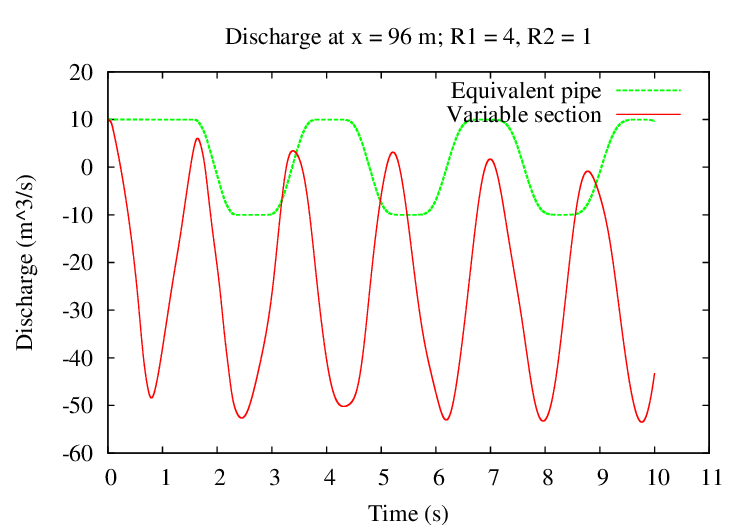} &
   \includegraphics[scale=0.59]{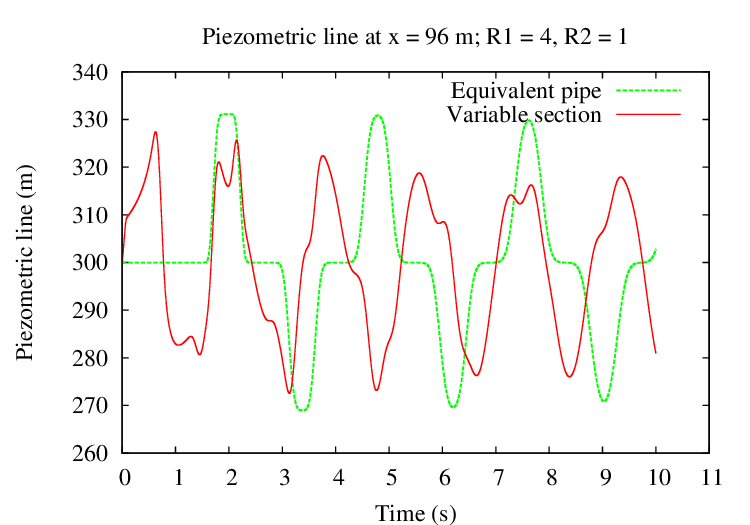} \\
 \end{tabular}
\end{center}
\caption{Discharge (left) and piezometric line (right) for $R_1 = 4\,m$}
\label{T1}
\end{figure}

\noindent We point out that the behaviour of the solutions corresponding to the equivalent pipe method and our method are different: this is due to the dynamic treatment of the term  $\displaystyle c^2 \frac{d \ln S}{dX}$ related to the variable section which is not present in the equivalent pipe method: see \fig{\ref{T12}}, \fig\ref{T11}, \fig\ref{T1} :

\begin{figure}[H]
\centering{
 \includegraphics[scale=0.9]{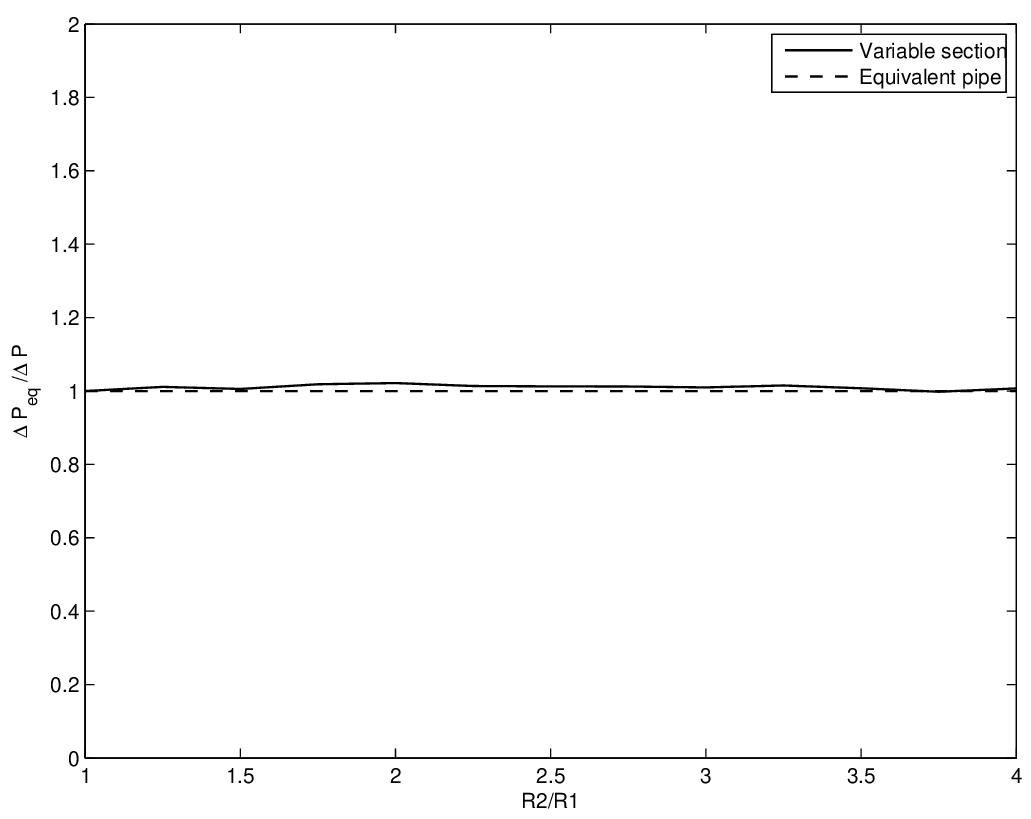}}
\caption{Comparison in the prediction of pressure rises in cone-shaped pipes between the present method and the equivalent pipe method}
\label{equivalent_pipe}
\end{figure}

%\nocite{*}
\bibliographystyle{plain}
%\bibliography{cinemix}

\end{document}